\documentclass[a4paper,12pt]{article}
\usepackage{amssymb,amsmath}
\def\beq{\begin{equation}}
\def\eeq{\end{equation}}
\def\bea{\begin{eqnarray}}
\def\eea{\end{eqnarray}}
\def\nn{\nonumber}
\def\U{{\cal U}}
\def\A{{\cal A}}
\def\R#1{{\cal R}_{#1}}

\def\hf{\frac{1}{2}}
\def\bra#1{\left\{ #1 \right\}_{q^{-1}}}
\def\ku#1{\left\{ #1 \right\}_q}
\def\qn#1{[#1]_q}
\def\sq#1{\big\lgroup #1 \big\rgroup_q}
\def\sqi#1{\big\lgroup #1 \big\rgroup_{q^{-1}}}
\def\uT#1{{\cal T}_{#1}}
\def\Tm#1#2{T^{#1}_{#2}}
\def\e#1#2{e^{#1}_{#2}}
\def\CGC#1#2{C^{#1}_{#2}\,}
\def\Pz#1#2{P^{#1}_{#2}(\zeta)}
\def\Fl#1#2#3{F^{\ell}(#1,#2,#3)\,}
\def\Gl#1#2#3{G^{\ell}(#1,#2,#3)\,}
\makeatletter
  \def\@cite#1#2{${\mbox{#1\if@tempswa , #2\fi}}$}
\makeatother
\def\ocite#1{$^{\mbox{\scriptsize{\cite{#1}}}}$}
\makeatletter
  \def\@biblabel#1{$^{\mbox{#1}}$}
\makeatother
\renewcommand{\theequation}{\arabic{section}.\arabic{equation}}

\begin{document}
%%%%%%%%%%%%%%%%%%%%%%%%%%%%%%%%%%%%%%%%%%%%%%%%%%%%%%%%%%%%%%
%
%  title page
%
%
%%%%%%%%%%%%%%%%%%%%%%%%%%%%%%%%%%%%%%%%%%%%%%%%%%%%%%%%%%%%%%
%
\thispagestyle{empty}

\vspace*{3cm}
\begin{center}
{\LARGE\sf
Universal ${\cal T}$-matrix, Representations of  \boldmath{$OSp_q(1/2)$} 
and Little \boldmath{$Q$}-Jacobi Polynomials}

\bigskip\bigskip
N. Aizawa

\bigskip
\textit{
Department of Mathematics and Information Sciences, \\
Graduate School of Science,\\
Osaka Prefecture University, \\
Daisen Campus, Sakai, Osaka 590-0035, Japan}\\
\bigskip
R. Chakrabarti, S.S. Naina Mohammed

\bigskip
\textit{
Department of Theoretical Physics, \\
University of Madras, \\
Guindy Campus, Chennai 600 025, India
}

\bigskip
J. Segar

\bigskip
\textit{
Department of Physics, \\
Ramakrishna Mission Vivekananda College, \\
Mylapore, Chennai 600 004, India
}
\bigskip

\end{center}

\vfill
\begin{abstract}
We obtain a  closed form expression of the universal ${\cal T}$-matrix 
encapsulating the duality between the quantum 
superalgebra $U_q[osp(1/2)]$ and the corresponding supergroup 
$OSp_q(1/2)$. The classical $q \rightarrow 1$ limit of this universal
${\cal T}$ matrix yields the group element of the undeformed
$OSp(1/2)$ supergroup. 
The finite dimensional representations of the quantum 
supergroup $OSp_q(1/2)$ are readily constructed employing 
the said universal ${\cal T}$-matrix and the known 
finite dimensional representations of the dually related deformed 
$U_q[osp(1/2)]$ superalgebra. Proceeding further, we derive 
the product law, the recurrence relations and the orthogonality of the 
representations of the quantum supergroup $OSp_q(1/2).$ It is shown that 
the entries of these representation matrices are expressed in terms 
of the {\it little $Q$-Jacobi polynomials with $Q = -q$}. 
Two mutually complementary 
singular maps of the universal ${\cal T}$-matrix on the universal 
${\cal R}$-matrix are also presented.
\end{abstract}
\newpage

\setcounter{page}{1}
%%%%%%%%%%%%%%%%%%%%%%%%%%%%%%%%%%%%%%%%%%%%%%%%%%%
%
%          Introduction 
%
%%%%%%%%%%%%%%%%%%%%%%%%%%%%%%%%%%%%%%%%%%%%%%%%%%%
%
\section{INTRODUCTION}
\label{Intro}

The representation theory of quantum groups and algebras have richer 
structures compared to their classical counterparts. 
Various non-classical features in the representation theory of the 
former have been found for specific values of deformation parameters 
such as roots of unity or crystal base limit. It is known that, 
for generic values of deformation parameters, each irreducible 
representation of the classical Lie groups and algebras has its quantum 
analog (see for example Ref.\cite{KS}). Even for such a generic case, 
however, there exist representations of the quantum algebra that do not 
have classical partners\ocite{MMNNSU,Ku91,NA93}. Extending our studies 
to the quantum supergroups, we expect further richness of 
representations as the nilpotency of Grassmann variables in 
classical supergroups are, in many cases, lost at quantum level. 
Grassmann coordinates of quantum superspaces\ocite{AC04} 
and quantum superspheres\ocite{AC05}, as well as the Grassmann 
elements of quantum supermatrices\ocite{AC04,AC05} are instances 
of lost nilpotency. When representations of such algebraic objects 
are considered, nonvanishing squares of Grassmann variables  
cause a drastic shift from the classical cases even for generic 
values of deformation parameters. 

  Influenced by this observation, we here study the representations of 
the simplest quantum supergroup $ OSp_q(1/2) $. 
Precise theory of matrix representations of quantum groups has been 
developed in Ref. \cite{Wor}.
Physical motivations of 
the present work is provided by the investigations on solvable models 
having quantized $osp(1/2)$ symmetry. For instance, vertex 
models\ocite{Sal}, Gaudin model\ocite{KM} and 2D field 
theories\ocite{SWK} based on $U_q[osp(1/2)]$ symmetry have been proposed 
and investigated. Fully developed representation theory of 
$OSp_q(1/2)$ will provide useful tools for analyzing these models and 
building new ones. 

  It is known that the representation theories of two quantum algebras 
$ U_q[sl(2)] $ and $ U_q[osp(1/2)] $ are quite parallel for a generic 
$q$. We naturally anticipate that the known results on the quantum group 
$SL_q(2)$ may be extended to the quantum supergroup $ OSp_q(1/2).$ 
To accomplish the extension, we employ the universal ${\cal T}$-matrix, 
which is a generalization of the exponential mapping relating a Lie 
algebra with its corresponding group. Underlying reasons for this are
as follows: (i) The universal ${\cal T}$ matrix succinctly embodies the 
representations of the dually related conjugate Hopf structures,
$ U_q[osp(1/2)] $ and $ OSp_q(1/2).$ In particular, contributions of 
the nonvanishing square of the odd elements of $ OSp_q(1/2)$ may be 
directly read from the expression of the universal ${\cal T}$-matrix. 
(ii) Moreover, the universal ${\cal T}$-matrix allows us to map each 
irreducible representation of $ U_q[osp(1/2)] $ on the corresponding 
one of $OSp_q(1/2). $ Therefore, various properties of representations 
follow from the corresponding ones of the universal ${\cal T}$-matrix, 
and the role of the lost nilpotency becomes explicit. Specifically, we 
demonstrate that the nonvanishing contributions of odd elements 
{\it assume polynomial structures} in the representation matrices. 

\par

The plan of this article is as follows. After fixing notations and 
conventions in the next section, the basis set of the Hopf dual to the 
$ U_q[osp(1/2)] $ algebra is explicitly obtained. The finitely generated 
basis sets of the dually related Hopf algebras are now used to derive  
a closed form expression of the universal ${\cal T}$-matrix via the 
method of Fr{\o}nsdal and Galindo\ocite{FG}. Singular, and, therefore, 
non-invertible maps of the universal ${\cal T}$-matrix on the universal 
${\cal R}$-matrix exist\ocite{FG}. Two such mutually complementary maps 
are studied in \S IV. General properties of the finite dimensional 
representations of the quantum supergroup $ OSp_q(1/2)$, such as the 
product law, the recurrence relations and the orthogonality of 
representations, follow, as observed in \S V, from the duality encompassed 
in the universal ${\cal T}$-matrix. Explicit form of the representation 
matrices are derived in \S VI and its relation to the little $Q$-Jacobi 
polynomials are discussed. It is shown that the entries of representation 
matrix are expressed in terms of the little $Q$-Jacobi polynomials with 
$Q = -q.$ This provides a new link of the representation 
theory of quantum supergroups with the hypergeometric functions. Section 
\S VII is devoted to concluding remarks. Corresponding results on 
$ SL_q(2)$ and other quantum groups are mentioned in each section.

%%%%%%%%%%%%%%%%%%%%%%%%%%%%%%%%%%%%%%%%%%%%%%%%%%%
%
%  U_q[osp(1/2)] and its representations
%
%%%%%%%%%%%%%%%%%%%%%%%%%%%%%%%%%%%%%%%%%%%%%%%%%%%
%
\setcounter{equation}{0}
\section{\boldmath{$U_q[osp(1/2)]$} AND ITS REPRESENTATIONS}
\label{UqRep}

  The quantum superalgebra $ \U = U_q[osp(1/2)] $ and the dually 
related quantum supergroup $ \A = OSp_q(1/2)$, dual to $\U$, have been 
introduced in Ref.\cite{KR}. Structures and representations of $ \U $  
have been investigated in Refs.\cite{Sal,KR}. For the purpose of fixing 
our notations and conventions we here list the relations that will be 
used subsequently. 

  The algebra $\U$ is generated by three elements $H$ (parity even) and 
$ V_{\pm} $ (parity odd) subject to the relations
\beq
   [H, V_{\pm}] = \pm\hf V_{\pm}, \qquad 
   \{ V_+, V_- \} = -\frac{q^{2H} - q^{-2H}}{q-q^{-1}} \equiv -\qn{2H}.
   \label{ospS2}
\eeq
The deformation parameter $q$ is assumed to be generic throughout this article. 
The Hopf algebra structures defined via the coproduct ($\Delta$), 
the counit ($\epsilon$) and the antipode ($S$) maps read as follows:
\bea
  & & \Delta(H) = H \otimes 1 + 1 \otimes H, \qquad
      \Delta(V_{\pm}) = V_{\pm} \otimes q^{-H} + q^H \otimes V_{\pm},
      \label{Dosp2} \\
  & & \epsilon(H) = \epsilon(V_{\pm}) = 0, \label{epS2} \\
  & & S(H) = -H, \qquad S(V_{\pm}) = -q^{\mp 1/2} V_{\pm}. \label{SospS2}
\eea
Using the flip operator $ \sigma $: 
$ \sigma(a \otimes b) = (-1)^{p(a)p(b)} b \otimes a$, where $p(a)$ 
denotes the parity of $ a, $ we define the transposed coproduct: 
$ \Delta' = \sigma \circ \Delta.$ The universal ${\cal R}$-matrix  
intertwining $ \Delta $ and $ \Delta'$ is given by\ocite{KR}
\beq
  \R{q} = 
  q^{4H \otimes H} \sum_{k \geq 0} 
  \frac{(q-q^{-1})^k q^{-k/2}}{\sqi{k}!} 
  (q^H V_+ \otimes q^{-H} V_-)^k,
  \label{UR2}  
\eeq
where
\beq
   \sq{k} = \frac{1-(-1)^k\, q^k}{1 + q},
   \quad
   \sq{k} ! = \sq{k} \sq{k-1} \cdots \sq{1},\quad
   \sq{0} ! = 1.
   \label{superqnumber}
\eeq
Two properties of $ \R{q}$ that will be used later are listed below: (i) 
It satisfies Yang-Baxter equation
\beq
   \R{q\,12} \R{q\,13} \R{q\,23} = \R{q\,23} \R{q\,13} \R{q\,12}, 
   \label{YBE}
\eeq
and (ii) its antipode map reads $ (S \otimes {\rm id}) \R{q} = 
\R{q^{-1}} = \R{q}^{-1}.$ 

  The finite dimensional irreducible representations of $ \U $ is specified by 
the highest weight $ \ell $ which takes any non-negative integral value. 
The irreducible representation space $ V^{(\ell)}$ of 
highest weight $\ell$ is $ 2\ell +1 $ dimensional. 
We denote its basis as 
$ \{\; e^{\ell}_m(\lambda) \ | \ m = \ell, \ell-1,\cdots, -\ell \;\}, $
where the index $ \lambda = 0, 1 $ specifies the parity of the highest weight 
vector $ e^{\ell}_{\ell}(\lambda). $ The parity of the vector 
$ e^{\ell}_m(\lambda) $ equals $ \ell - m + \lambda,$ 
as it is obtained by the action of $ V_-^{\ell-m} $ on 
$ e^{\ell}_{\ell}(\lambda).$  For the superalgebras the norm of the 
representation basis need not be chosen positive definite.  
In this work, however, we assume the positive definiteness of 
the basis elements: 
\beq
 (e^{\ell}_m(\lambda), e^{\ell'}_{m'}(\lambda)) = \delta_{\ell \ell'} 
 \delta_{mm'}.
 \label{Normalization}
\eeq
With these settings, the irreducible representation of $ \U$ on $ V^{(\ell)} $ 
is given by
\bea
  & & H \e{\ell}{m}(\lambda) = \frac{m}{2} \e{\ell}{m}(\lambda), \nn \\
  & & V_+ \e{\ell}{m}(\lambda) = \left(\frac{1}{\ku{2}}\ku{\ell-m}\ku{\ell+m+1}\right)^{1/2} 
      \e{\ell}{m+1}(\lambda),
      \label{RepU} \\
  & & V_- \e{\ell}{m}(\lambda) 
     = (-1)^{\ell-m-1} \left(\frac{1}{\ku{2}} \ku{\ell+m}\ku{\ell-m+1}\right)^{1/2}
      \e{\ell}{m-1}(\lambda), \nn
\eea
where
\beq
  \ku{m} = \frac{q^{-m/2}-(-1)^m q^{m/2}}{q^{-1/2}+q^{1/2}}.  \label{Ksymb}
\eeq

  Tensor product of two irreducible representations is, in general, 
reducible and may be decomposed into irreducible ones without 
multiplicity:
\[
  V^{(\ell_1)} \otimes V^{(\ell_2)} = 
  V^{(\ell_1+\ell_2)} \oplus V^{(\ell_1+\ell_2 -1)} \oplus \cdots 
  \oplus V^{(|\ell_1-\ell_2|)}.
\]
The decomposition of the tensored vector space in the irreducible basis 
is provided by the Clebsch-Gordan coefficients (CGC):
\beq
    \e{\ell}{m}(\ell_1,\ell_2,\Lambda) = \sum_{m_1,m_2} 
    \CGC{\ell_1\ \ell_2\ \;\ell}{m_1\,m_2\,m}\,\,
    \e{\ell_1}{m_1}(\lambda) \otimes \e{\ell_2}{m_2}(\lambda),
    \label{CGCdef}
\eeq
where $m= m_1+m_2, $ and 
$ \Lambda = \ell_1+\ell_2+\ell \ ({\rm mod}\ 2) $ signify the parity of 
the highest weight vector $ \e{\ell}{\ell}(\ell_1,\ell_2,\Lambda).$ 
The CGC for $\U$ is extensively studied in Ref.\cite{MM}. In spite of 
our assumption (\ref{Normalization}) regarding the positivity of 
the basis vectors, the norm of tensored vector space is not 
always positive definite. 
Indeed, the basis (\ref{CGCdef}) is pseudo orthogonal:
\beq
  (\e{\ell'}{m'}(\ell_1,\ell_2,\Lambda), \e{\ell}{m}(\ell_1,\ell_2,\Lambda))
  = 
  (-1)^{(\ell-m+\lambda)(\ell_1+\ell_2+\ell+\lambda)}
  \delta_{\ell'\ell}\delta_{m'm}.
  \label{Normprod}
\eeq
The CGC satisfies two pseudo orthogonality relations
\beq
  \sum_{m_1,m_2} (-1)^{(\ell_1-m_1+\lambda)(\ell_2-m_2+\lambda)} 
  \CGC{\ell_1\ \ell_2\ \;\ell}{m_1\,m_2\,m} 
  \CGC{\ell_1\ \ell_2\ \;\ell'}{m_1\,m_2\,m'} 
  =
  (-1)^{(\ell-m+\lambda)(\ell_1+\ell_2+\ell+\lambda)} \delta_{\ell\ell'}\delta_{mm'},
  \label{CGCortho1}
\eeq
\beq
  \sum_{\ell,m}(-1)^{(\ell-m+\lambda)(\ell_1+\ell_2+\ell+\lambda)}
  \CGC{\ell_1\ \ell_2\ \;\ell}{m_1\,m_2\,m} 
  \CGC{\ell_1\ \ell_2\ \; \ell}{m_1'\,m_2'\,m}
  =
  (-1)^{(\ell_1-m_1+\lambda)(\ell_2-m_2+\lambda)} 
  \delta_{m_1\,m_1'} \delta_{m_2\,m_2'}.
  \label{CGCortho2}
\eeq
Equation (\ref{CGCortho2}) immediately provides  the 
inversion of  the construction (\ref{CGCdef}): 
\beq
   \e{\ell_1}{m_1}(\lambda) \otimes \e{\ell_2}{m_2}(\lambda) 
   = 
   (-1)^{(\ell_1-m_1)(\ell_2-m_2)} \sum_{\ell,m} (-1)^{(\ell-m)(\ell_1+\ell_2+\ell)}
   \CGC{\ell_1\ \ell_2\ \;\ell}{m_1\,m_2\,m} 
   \,\,\e{\ell}{m}(\ell_1,\ell_2,\Lambda).
   \label{CGCreverse}
\eeq
Before closing this section, we make two remarks: (i) All the CGC are of 
parity even state. (ii) The explicit realization of CGC for $\U$ is found 
in Ref.\cite{AC05}, and also in Ref.\cite{MM}. As we maintain the phase 
convention for the representation of $ \U $ given in Ref.\cite{MM}, we 
use the results obtained therein. 
 
%%%%%%%%%%%%%%%%%%%%%%%%%%%%%%%%%%%%%%%%%%%%%%%%%%
%
%    Dual algebra and universal T-matrix
%
%%%%%%%%%%%%%%%%%%%%%%%%%%%%%%%%%%%%%%%%%%%%%%%%%%
%
\setcounter{equation}{0}
\section{UNIVERSAL $\cal T$-MATRIX VIA DUALITY}
\label{DualT}
Two Hopf algebras $\U$ and $\A$ are in duality \ocite{FG} if there 
exists a doubly-nondegenerate bilinear form
$
\langle \; ,\; \rangle: \A \otimes \U 
\ \rightarrow \  {\mathbb C} 
% \qquad \forall {\sf a} \in \A,\,\,\forall {\sf u} \in \U,
%\label{au_form}
$
such that, for $({\sf a}, {\sf b}) \in \A, 
({\sf u}, {\sf v}) \in \U$,
\bea
&&\langle {\sf a}, {\sf u v}\rangle = \langle \Delta_{\A}({\sf a}),
{\sf u} \otimes {\sf v}\rangle, \qquad
\langle {\sf a b}, {\sf u}\rangle = \langle {\sf a}\otimes {\sf b},
\Delta_{\U}({\sf u})\rangle,\nn\\
&&\langle {\sf a}, 1_{\U}\rangle = \epsilon_{\A}({\sf a}),\quad 
\langle 1_{\A}, {\sf u}\rangle = \epsilon_{\U}({\sf u}),\quad 
\langle {\sf a}, S_{\U}({\sf u})\rangle =  
\langle S_{\A}({\sf a}), {\sf u}\rangle.   
\label{dual_def}
\eea
Let the ordered monomials $E_{k \ell m} = V_{+}^{k} 
H^{\ell} V_{-}^{m},\,\, (k,\ell,m) \in (0, 1, 2,\cdots)$ be the 
basis elements of the algebra $\U$ obeying the multiplication and 
the induced coproduct rules given by 
\beq
  E_{k\ell m} E_{k'\ell'm'} = \sum_{pqr} f^{\quad pqr}_{k\ell m \ k'\ell'm'} E_{pqr},
  \qquad
  \Delta(E_{k\ell m}) = \sum_{pqr \atop p'q'r'} g^{pqr\ p'q'r'}_{\ k\ell m} 
  E_{pqr} \otimes E_{p'q'r'}.
  \label{strE}
\eeq
The basis elements $e^{k \ell m}$ of the dual Hopf algebra $\A$  
follow the relation
\beq
  \langle e^{k \ell m}, E_{k'\ell'm'} \rangle = 
  \delta^k_{k'} \delta^{\ell}_{\ell'} \delta^m_{m'}.
  \label{pair}
\eeq
In particular, the generating elements of the algebra $\A$, defined as 
$x = e^{1 0 0}, y = e^{0 0 1}$ and $z = e^{0 1 0}$, satisfy the 
following duality structure:
\beq
\langle x, V_+\rangle = 1, \quad
\langle z, H\rangle = 1, \quad
\langle y, V_-\rangle = 1. 
\label{xyz_def}
\eeq
Thus, $x$ and $y$ are of odd parity, while $z$ is even. 
The duality condition (\ref{dual_def}) requires the basis set $e^{k \ell m}$ 
to obey the multiplication and coproduct rules given below: 
\beq
   e^{pqr} e^{p'q'r'} = \sum_{k\ell m} g^{pqr\ p'q'r'}_{\quad k\ell m} e^{k\ell m},
   \qquad
   \Delta(e^{pqr}) = \sum_{k\ell m \atop k'\ell'm'} f^{\quad pqr}_{k\ell m\ k'\ell'm'} 
   e^{k\ell m} \otimes e^{k'\ell'm'}.
   \label{stre}
\eeq

  To derive the Hopf properties of the dual algebra $\A$, we, therefore, 
need to extract the structure constants defined in (\ref{strE}).
Towards this end we note that the induced coproduct map of the 
elements $E_{k \ell m}$ may be obtained via (\ref{Dosp2}):
\bea
  \Delta(E_{k \ell m}) &=& \Delta(V_+)^{k} \Delta(H)^{\ell} \Delta(V_-)^m 
  \nn \\
  &=& \sum_{a=0}^{k} \sum_{b=0}^{\ell} \sum_{c=0}^{m} 
  \big\lgroup \begin{array}{c} k \\ a \end{array} \big\rgroup_q
  \left( \begin{array}{c} \ell \\ b \end{array} \right)
  \big\lgroup \begin{array}{c} m \\ c \end{array} \big\rgroup_q
  (-1)^{(m-c)(a+c)} q^{-a(k-a)/2-c(m-c)/2}
  \nn \\
  &\times& V_+^{k-a} q^{(a+c)H} H^{\ell-b} V_-^{m-c} 
      \otimes
      V_+^a q^{-(k+m-a-c)H} H^b V_-^c,
  \label{ExpandDE}
\eea
where
\[
  \big\lgroup  \begin{array}{c} k \\ a \end{array} \big\rgroup_q
  = \frac{ \sq{k}! }{ \sq{a}! \sq{k-a}! }.
\]
The second equality in (\ref{ExpandDE}) can be verified by using the 
commutation relations
\bea
  & &  q\, (V_+ \otimes q^{-H}) (q^H \otimes V_+) + (q^H \otimes V_+) (V_+ \otimes q^{-H}) = 0, \nn \\
  & &  q\, (q^H \otimes V_-) (V_- \otimes q^{-H}) + (V_- \otimes q^{-H}) (q^H \otimes V_-) = 0, \nn
\eea
and an extension of binomial theorem that is easily proved by 
induction: Arbitrary operators $ A, B $ subject to the commutation 
properties $ q AB + BA = 0,$ satisfy the following expansion 
\beq
  (A+B)^n = \sum_{k=0}^n \big\lgroup \begin{array}{c} n \\ k \end{array} \big\rgroup_q
  A^{n-k} B^k. 
  \label{binomial}
\eeq
Employing ({\ref{ExpandDE}) we now obtain a set of structure constants: 
\bea
  & & g^{100\ 001}_{\ \ k\ell m} = \delta_{k 1} \delta_{\ell 0} \delta_{m 1}, 
  \qquad
  g^{001\ 100}_{\ \ k\ell m} = -\delta_{k 1} \delta_{\ell 0} \delta_{m 1}, 
  \nn \\
  & & g^{100\ 010}_{\ \ k\ell m} = -\ln q \, \delta_{k 1} \delta_{\ell 0} \delta_{m 0}
  + \delta_{k1} \delta_{\ell 1} \delta_{m 0},
  \nn \\
  & &   g^{010\ 100}_{\ \ k\ell m} = \ln q \, \delta_{k 1} \delta_{\ell 0} \delta_{m 0}
  + \delta_{k1} \delta_{\ell 1} \delta_{m 0},
  \nn \\
  & & g^{010\ 001}_{\ \ k\ell m} = \ln q \, \delta_{k 0} \delta_{\ell 0} \delta_{m 1}
  + \delta_{k 0} \delta_{\ell 1} \delta_{m 1},
  \nn \\
  & &   g^{001\ 010}_{\ \ k\ell m} = -\ln q \, \delta_{k 0} \delta_{\ell 0} \delta_{m 1}
  + \delta_{k 0} \delta_{\ell 1} \delta_{m 1}.
  \nn 
\eea
The above structure constants immediately yield the algebraic 
relations obeyed by the generators of the algebra $\A$:
\beq
  \{ x, y \} = 0, \qquad [z, x] = 2 \ln q\,\,x, 
  \qquad [z, y] = 2 \ln q\,\, y.
  \label{Acomm}
\eeq

  Proceeding towards constructing the coproduct maps of the 
generating elements of the dual algebra $\A$ we notice that the defining 
properties (\ref{stre}) provide the necessary recipe:
\bea
\Delta(x)&=& \sum_{k \ell m \, \atop k^{\prime} \ell^{\prime} m^{\prime}}
\, f^{\ \;1 0 0}_{k \ell m\; k^{\prime} \ell^{\prime} m^{\prime}} \, 
e^{k \ell m} \, \otimes 
e^{k^{\prime} \ell^{\prime} m^{\prime}},\nn\\
\Delta(z)&=& \sum_{k \ell m \, \atop k^{\prime} \ell^{\prime} m^{\prime}}\, 
f^{\ \;0 1 0}_{k \ell m\; k^{\prime} \ell^{\prime} m^{\prime}}  \,  
e^{k \ell m} \, \otimes e^{k^{\prime} \ell^{\prime} m^{\prime}}, \nn\\
\Delta(y)&=& \sum_{k \ell m \, \atop k^{\prime} \ell^{\prime} m^{\prime}}\, 
f^{\ \;0 0 1}_{k \ell m\; k^{\prime} \ell^{\prime} m^{\prime}} 
\, e^{k \ell m} \, \otimes e^{k^{\prime} \ell^{\prime} m^{\prime}}. 
\label{xyz_copro}
\eea
The relevant structure constants obtained via (\ref{strE}) are listed 
below: 
\bea
  & & f^{\ \;1 0 0}_{k \ell m\; k^{\prime} \ell^{\prime} m^{\prime}}
      = \delta_{k1}\delta_{\ell 0}\delta_{m0} \delta_{k'0}\delta_{\ell' 0}\delta_{m'0}
      + (-1)^m \frac{1}{2^{\ell}} \sigma_{m+1} \delta_{k0}\delta_{k'\; m+1}\delta_{\ell'0}\delta_{m'0},
  \nn \\
  & & f^{\ \;0 1 0}_{k \ell m\; k^{\prime} \ell^{\prime} m^{\prime}} 
      = \delta_{k0}\delta_{m0}\delta_{k'0}\delta_{m'0} (\delta_{\ell 1}\delta_{\ell'0}
         + \delta_{\ell 0}\delta_{\ell'1})
         + (-1)^m \frac{4\ln q}{q-q^{-1}} \sigma_m 
         \delta_{k0}\delta_{\ell 0} \delta_{k'm} \delta_{\ell'0}\delta_{m'0},
   \nn \\
   & & f^{\ \;0 0 1}_{k \ell m\; k^{\prime} \ell^{\prime} m^{\prime}} 
      = \delta_{k0}\delta_{\ell 0}\delta_{m0} \delta_{k'0}\delta_{\ell'0}\delta_{m'1} 
      + (-1)^{k'}\frac{1}{2^{\ell'}} \sigma_{k'+1} \delta_{k0}\delta_{\ell0}\delta_{m\; k'+1} 
        \delta_{m'0},
   \label{SCdualA} \\
   & & \sigma_1 = 1, \qquad \sigma_m = \prod_{k=1}^{m-1} \sum_{\ell=0}^{k-1} (-1)^{\ell}\,
       \qn{k-\ell}, \quad (m >1).
   \nn
\eea
The coproduct maps of the generators of $\A$ may now be explicitly 
obtained {\it {\`a} la} (\ref{xyz_copro}) provided the basis elements
$e^{k \ell m}$ of the  algebra $\A$ are known. We complete this 
task subsequently. 

  As the algebra $\A$ is finitely generated, we may start with the 
generators $(x, y, z)$ and obtain all dual basis elements
$e^{k \ell m},\;(k, \ell, m) \in (0, 1, 2,\cdots)$ by successively 
applying the multiplication rule given in the first equation in 
(\ref{stre}). The necessary structure constants may be read 
from the relation (\ref{strE}) of the algebra $\U$. In the procedure 
described below we maintain the operator ordering of the monomials as 
$x^{k} z^{\ell} y^{m},\,\, (k, \ell, m) \in (0, 1, 2,\cdots)$. 
The product rule 
\beq
e^{1 0 0} \,\, e^{n 0 0} = \sum_{k \ell m}\, 
g^{1 0 0\,n 0 0}_{\ \; k \ell m}\, e^{k \ell m}
\label{prod_1} 
\eeq
and the explicit evaluation of the structure constant
\beq
  g^{100 \ n00}_{\ \ k\ell m}
  = \ku{n+1} \delta_{k\; n+1} \delta_{\ell 0} \delta_{m 0} 
  \label{strC-x}
\eeq
obtained from (\ref {ExpandDE}) immediately provide
\beq
   e^{n00} = \frac{x^n}{\ku{n}!}, \qquad
   \ku{n}! = \prod_{\ell=1}^{n} \ku{\ell}, \quad \ku{0}! = 1. \label{en00}
\eeq
Employing another product rule
\beq
e^{n r 0} \, e^{0 1 0} = \sum_{k \ell m} g^{n r 0 \; 0 1 0}_{\ \;k \ell m} 
\,\, e^{k \ell m}
\label{prod_2}
\eeq   
and the value of the relevant structure constant 
\[
  g^{nr0\ 010}_{\ k\ell m} = -n \ln q \,\delta_{kn} \delta_{\ell r} \delta_{m 0} 
  + (r+1) \delta_{k n} \delta_{\ell\; r+1} \delta_{m 0},
\]
obtained in the aforesaid way we produce the following result:
\beq
  e^{nr0} = \frac{x^n}{\ku{n}!} \frac{1}{r!} (z+ n \ln q)^r.   \label{enr0}
\eeq
Continuing the above process of building of the dual basis set we use 
the product rule
\beq
e^{n r s} \, e^{0 0 1} = \sum_{k \ell m} g^{n r s \, 0 0 1}_{\ \;k \ell m} 
\,\, e^{k \ell m}
\label{prod_3}
\eeq   
and the value of the corresponding structure constant
\[
  g^{nrs\ 001}_{\ \ k\ell m} = \bra{s+1} \sum_{j=0}^r \frac{1}{j!}\, (\ln q)^j 
  \delta_{kn} \delta_{\ell\;r-j} \delta_{m\;s+1},
\]
obtained via (\ref{ExpandDE}). This finally leads us to the complete 
construction of the basis element of the algebra $\A$:
\beq
   e^{nrs} = \frac{x^n}{\ku{n}!} \frac{(z+(n-s) \ln q)^r}{r!} \frac{ y^s}{\bra{s}!}.
   \label{enrs}
\eeq
Combining our results in (\ref{xyz_copro}), (\ref{SCdualA}) and 
(\ref{enrs}), we now provide the promised coproduct structure of the 
generators of the algebra $\A$:
\bea
  & &   \Delta(x) = x \otimes 1 + \sum_{m=0}^{\infty} (-1)^m \sigma_{m+1} 
  q^{-m/2} e^{z/2} \frac{y^m}{\bra{m}!} \otimes \frac{x^{m+1}}{\ku{m+1}!},
  \nn \\
  & &   \Delta(z) = z \otimes 1 + 1 \otimes z + 
  \frac{4 \ln q}{q-q^{-1}} \sum_{m=1}^{\infty} (-1)^m \sigma_m 
  \frac{y^m}{\bra{m}!} \otimes \frac{x^m}{\ku{m}!},
  \label{Deltaz} \\
  & &    \Delta(y) = 1 \otimes y + \sum_{k=0}^{\infty} (-1)^k \sigma_{k+1} 
   q^{k/2} \frac{y^{k+1}}{\bra{k+1}!} \otimes \frac{x^k}{\ku{k}!} e^{z/2}.
  \nn
\eea
Algebraic simplifications allow us to express the 
coproduct maps of the above generators more succinctly:
\bea
  & & \Delta(x) = x \otimes 1 + 
      \sum_{m=0}^{\infty} (-1)^{m(m-1)/2}
      \left(
         \frac{1+q^{-1}}{q-q^{-1}}
      \right)^m  e^{z/2} y^m \otimes x^{m+1},
      \nn \\
  & & \Delta(z) = z \otimes 1 + 1 \otimes z 
      + \frac{4 \ln q}{q-q^{-1}} \sum_{m=1}^{\infty} 
      \frac{(-1)^{m(m+1)/2}}{\bra{m}}
      \left(
         \frac{q^{1/2}+q^{-1/2}}{q-q^{-1}}
      \right)^{m-1} y^m \otimes x^{m},
      \nn \\
  & & \label{Deltaz2} \\
  & & \Delta(y) = 1 \otimes y + \sum_{m=0}^{\infty} (-1)^{m(m-1)/2} 
      \left(
         \frac{q+1}{q-q^{-1}}
      \right)^m  y^{m+1} \otimes x^{m} e^{z/2}.
      \nn
\eea
With the aid of the result (\ref{Deltaz2}) we may explicitly 
demonstrate that the coproduct map is a homomorphism of the algebra
(\ref{Acomm}): namely,
\[
  \{ \Delta(x), \Delta(y) \} = 0, \qquad
  [\Delta(z), \Delta(x)] = 2 \ln q\, \Delta(x), \qquad
  [\Delta(z), \Delta(y)] = 2 \ln q\, \Delta(y).  
\]
The coassociativity constraint 
\[
(\hbox{id} \otimes \Delta) \circ \Delta({\cal X}) =
(\Delta \otimes \hbox{id}) \circ \Delta({\cal X}) \qquad
\forall{\cal X} \in (x, y, z)
\]
may also be established by using the following identity:
\bea
&& \exp (\Delta (z)) = (\exp(z) \otimes 1)\,\,
\prod_{m = 1}^{\infty}\,{\cal P}_{m}\,\,
(1 \otimes \exp (z)),\nn\\
&&{\cal P}_{m} = \exp \left((- 1)^{m (m + 1)/2}\,\,\frac{[2 m]_{q}}
{m\,\{m\}_{q^{-1}}}\,\,\Big(\frac{q^{1/2} + q^{- 1/2}}{q - q^{-1}}
\Big)^{m -1}
\,\,y^{m}\, \otimes \, x^{m}\right).
\label{exp_del_z}
\eea   
 The counit map of the generators of the algebra $\A$ reads as
\beq
\epsilon (x) = \epsilon (y) = \epsilon (z) = 0.
\label{xyz_counit}
\eeq
The antipode map of the dual generators follows from the 
last equation in (\ref{dual_def}). We quote the results here: 
\bea
 & & S(x) = -\sum_{m=0}^{\infty} (-1)^{m(m-1)/2} q^{-1} 
     \left(
        \frac{1+q^{-1}}{q-q^{-1}}
     \right)^m x^{m+1} \exp\left( -\frac{m+1}{2}z \right) y^m,
     \nn \\
 & & S(z) = -z + \frac{4 \ln q}{q-q^{-1}} \sum_{m=1}^{\infty} (-1)^{m(m+1)/2} 
     \frac{1}{\bra{m}} 
     \left(
         \frac{q^{1/2} + q^{-1/2}}{q-q^{-1}}
     \right)^{m-1} x^m e^{-mz/2} y^m,
     \nn \\ %\label{Sz2} 
  & & \label{Sz2} \\
  & & S(y) = -\sum_{m=0}^{\infty} (-1)^{m(m+1)/2} q 
     \left(
       \frac{q+1}{q-q^{-1}}
     \right)^m x^m \exp\left( -\frac{m+1}{2} z \right) y^{m+1}.
     \nn
\eea
This completes our construction of 
the Hopf algebra $\A$ dually related to the quantum superalgebra 
$\U$.

  Our explicit listing of the complete set of dual basis elements in 
(\ref{enrs}) allows us to obtain {\it \`{a} la} Fr\o nsdal and 
Galindo \ocite {FG} the universal ${\cal T}$-matrix for the supergroup:
\beq
   {\cal T}_{e,E} = 
   \sum_{k\ell m} (-1)^{p(e^{k\ell m}) (p(e^{k\ell m})-1)/2}\, 
   e^{k\ell m} \otimes E_{k\ell m},
   \label{uniTdef}
\eeq
where the parity of basis elements is same for two Hopf algebras $\U $ and $\A$
\beq
   p(e^{k\ell m}) = p(E_{k\ell m}) = k + m.   \label{paritygen}
\eeq
The notion of the universal ${\cal T}$-matrix is a key feature capping 
the Hopf duality structure. Consequently, the duality relations 
(\ref{dual_def}) may be concisely expressed \ocite{FG} in terms of the 
${\cal T}$-matrix as
\bea
    & & \uT{e,E} \uT{e',E} = \uT{\Delta(e),E}, \qquad 
        \uT{e,E} \uT{e,E'} = \uT{e,\Delta(E)}, \nn \\
    & & \uT{\epsilon(e),E} = \uT{e,\epsilon(E)} = 1, \quad \ \;
        \uT{S(e),E} = \uT{e,S(E)}.  \label{uniTprop}
\eea
where $e$ and $e^{\prime}$\, ($E$ and $E^{\prime}$) refer to the two 
identical copies of algebra $\A$ ($\U$). 
A general discussion of the universal ${\cal T}$-matrix for supergroups 
is found in the Appendix.

  As both the Hopf algebras in our case are finitely generated, the 
universal ${\cal T}$-matrix may now be obtained as an operator valued 
function in a closed form: 
\bea
  \uT{e,E} &=& 
  \left( \sum_{k=0}^{\infty} \frac{(x \otimes V_+q^H)^k}{\sq{k}!} \right) 
  \exp(z\otimes H) 
  \left(
     \sum_{m=0}^{\infty} \frac{(y \otimes q^{-H} V_-)^m}{\sqi{m}!} 
  \right)
  \nn \\
  &\equiv& {}_{\times}^{\times} 
  {\cal E}{\rm xp}_q (x \otimes V_+q^H) \; \exp(z \otimes H) \;
  {\cal E}{\rm xp}_{q^{-1}} (y \otimes q^{-H}V_-){}_{\times}^{\times},
  \label{uniTclosed}
\eea
where we have introduced a deformed exponential that is characteristic 
of the quantum $OSp_{q}(1/2)$ supergroup:
\beq
  {\cal E}{\rm xp}_q(x) \equiv \sum_{n=0}^{\infty} \frac{x^n}{\sq{n}!},
  \label{deformedExp}
\eeq
The operator ordering has been explicitly indicated in (\ref{uniTclosed}). 
The closed form of the universal ${\cal T}$-matrix in (\ref{uniTclosed}) 
will be used in the computation of representation matrices of 
the quantum supergroup $\A$. In Ref. \cite{DKLS}, using the Gauss
decomposition of the fundamental representation a universal 
$\cal T$-matrix for $ \U $ is given in terms of the 
standard $q$-exponential instead of the deformed exponential 
(\ref{deformedExp}) characterizing quantum supergroups.

\par 

  In the classical limit of $ q \rightarrow 1,$ it is immediately 
evident that the structure constant (\ref{strC-x}) is truncated at 
$ n=2 $ so that $x$ remains nilpotent. Similarly, $y^2 = 0 $ holds in 
this limit. It is interesting to observe the $q \rightarrow 1$ limit of 
the universal ${\cal T}$ matrix (\ref{uniTclosed}). For this 
purpose we note   
\beq
\lim_{q \rightarrow 1} \sq{2 n} \rightarrow \,n (1 - q),\qquad
\lim_{q \rightarrow 1} \sq{2 n + 1} \rightarrow 1, \quad 
n = 0, 1, 2,\cdots.
\label{sqlim}
\eeq
Assuming the finite limit 
\beq
\lim_{q \rightarrow 1} \frac{x^{2}}{q - 1} = \mathfrak{x},\qquad
\lim_{q \rightarrow 1} \frac{y^{2}}{q^{- 1} - 1} = \mathfrak{y},
\label{xylim}
\eeq
it immediately follows that the universal ${\cal T}$ matrix 
(\ref{uniTclosed}) reduces to the group element of the undeformed
$OSp(1/2)$, and by definition constitutes its universal 
${\cal T}$-matrix: 
\beq
 \uT{e,E} \big|_{q \rightarrow 1} = (1 \otimes 1 + x \otimes V_+ )\, 
 \exp (\mathfrak{x} \otimes V_{+}^{2}) \,\exp(z \otimes H) 
 \exp (\mathfrak{y} \otimes V_{-}^{2}) \, 
 (1\otimes 1 + y \otimes V_-).
 \label{uniTclassical}
\eeq
The finite limiting elements $(\mathfrak{x, y})$ are bosonic in nature, 
and in the context of the classical limit of the function algebra 
${\cal A}$ they are dually related to the squares of the odd generators 
of the undeformed $osp(1/2)$. The elements $(V_{\pm}^{2}, H)$ of the 
classical $osp(1/2)$ algebra form a $sl(2)$ subalgebra. The corresponding 
classical $SL(2)$ subgroup structure is evident from 
(\ref{uniTclassical}). In 
fact, the correct limiting structure (\ref{uniTclassical}) emphasizes
the essential validity of the quantum universal ${\cal T}$ matrix
derived in (\ref{uniTclosed}).    
Obviously, there is a striking difference between the 
quantum and classical universal ${\cal T}$-matrices caused by the 
absence of nilpotency of parity odd elements in the former case. 
The infinite series of operators summarized in the deformed 
exponential contribute to new polynomial matrix elements in the 
representations of quantum supergroup $\A$. This is  considered  
in detail in \S V. 

  The above construction of the universal ${\cal T}$-matrix for the 
algebra $\U$ is parallel to the one for the generalized Heisenberg 
algebra\ocite{ACS1} which is a bosonization\ocite{McMaj} of the 
superalgebra $\U.$ Dual basis to the two-parametric deformation of 
$ GL(2) $ is studied in Ref. \cite{Dov}. The universal 
${\cal T}$-matrix for the two-parametric quantum $ GL(2) $ is given 
in Refs.\cite{FG,JV1}. The generalization to the quantum $ gl(n) $ is 
found in Ref.\cite{Fro} and 
a supersymmetric extension is initiated in Ref.\cite{CJ1}.

%%%%%%%%%%%%%%%%%%%%%%%%%%%%%%%%%%%%%%%%%%%%%%%%%%
%
%   Mappings to universal R-matrix
%
%%%%%%%%%%%%%%%%%%%%%%%%%%%%%%%%%%%%%%%%%%%%%%%%%%
%
\setcounter{equation}{0}
\section{MAPPING ${\cal T}$ ON ${\cal R}$}
\label{T2R}

  Two singular and mutually complementary maps connecting the universal 
${\cal T}$-matrix in (\ref{uniTclosed}) and the universal 
${\cal R}$-matrix in (\ref{UR2}) are discussed in this section. The 
first map $ \Phi: {\cal A} \rightarrow {\cal U} $ reads
\beq
  \Phi(x) = 0, \qquad 
  \Phi(z) = (4 \ln q)H, \qquad
  \Phi(y) = q^{-1/2} (q-q^{-1}) q^H V_+. \label{map1}
\eeq 
It is easily observed to satisfy the following properties: 
(i) $ (\Phi \otimes {\rm id})(\uT{e,E}) = \R{q},$ and
(ii) $\Phi$ is an algebra homomorphism though not a Hopf 
algebra homomorphism; that is, $ \Phi $ respects the 
commutation relations in (\ref{Acomm}) but does not maintain
the Hopf coalgebra maps. 

  For introducing the second map, we recast the 
universal ${\cal T}$-matrix in the form given below:
\bea
  \uT{E,e} 
  &=& \sum_{k\ell m} (-1)^{(k+m)(k+m-1)/2} E_{k\ell m} \otimes e^{k\ell m}
  \nn \\
  &=& \left(
     \sum_k \frac{1}{\sq{k}!} (V_+q^H \otimes x)^k
  \right) 
  e^{H \otimes z} 
  \left(
      \sum_m \frac{1}{\sqi{m}!} (q^{-H} V_- \otimes y)^m
  \right).
  \label{uniTEe}
\eea
The universal ${\cal R}$-matrix is also rewritten as
\beq
  \R{q} = 
   \sum_{k \geq 0} 
  \frac{(q-q^{-1})^k }{\sqi{k}!} 
  (V_+ q^{-H} \otimes q^{H} V_-)^k q^{4H \otimes H}.
  \label{UR3}  
\eeq
Comparison of above two expressions immediately yields the promised map 
$ \Psi : {\cal A} \rightarrow {\cal U} $ defined as follows:
\beq
  \Psi(x) = (q^{-1}-q)\, q^{-H} V_-, \qquad \Psi(z) = (-4\ln q) H, 
  \qquad \Psi(y) = 0. \label{map2}
\eeq
One can immediately verify that 
$ ({\rm id} \otimes \Psi) (\uT{E,e}) = \R{q}^{-1},$ 
and that $ \Psi $ is an algebraic homomorphism but not a 
Hopf algebra homomorphism. It is interesting to observe that in both the 
maps introduced here {\it one} Borel subalgebra of the function algebra 
${\cal A}$ is mapped on the corresponding Borel subalgebra of the 
universal enveloping algebra ${\cal U}$. Therefore, the two conjugate 
Borel subalgebras of the ${\cal U}$ algebra are acted upon by two 
distinct, but complementary maps. Being singular in nature, these maps 
are, however, not invertible.    

  The maps $ \Phi $ and $\Psi $ may be utilized to 
connect the universal ${\cal T}$-matrix and the Yang-Baxter equation. 
As shown in the Appendix, the universal ${\cal T}$-matrix 
satisfies $RTT$-type relations. Using the tensored operators
\bea
  & & \uT{e,E}^{(1)} = \sum_k (-1)^{(k+m)(k+m-1)/2} 
  e^{k\ell m} \otimes E_{k\ell m} \otimes 1,
  \nn \\
  & & \uT{e,E}^{(2)} = \sum_k (-1)^{(k+m)(k+m-1)/2} 
  e^{k\ell m} \otimes 1 \otimes E_{k\ell m}, \label{Ttensor1}
\eea  
the following identity may be established:
\beq
  (1 \otimes \R{q})\, \uT{e,E}^{(1)}\, \uT{e,E}^{(2)}
  =
  \uT{e,E}^{(2)}\, \uT{e,E}^{(1)} \, (1 \otimes \R{q}).
  \label{RTT1}
\eeq
Mirroring the structure in (\ref{Ttensor1}) we also define the 
transposed ${\cal T}$-matrices in the tensored space as 
\bea
  & & \uT{E,e}^{(1)} = \sum_k (-1)^{(k+m)(k+m-1)/2} 
  E_{k\ell m} \otimes 1 \otimes e^{k\ell m},
  \nn \\
  & & \uT{E,e}^{(2)} = \sum_k (-1)^{(k+m)(k+m-1)/2} 
  1 \otimes E_{k\ell m} \otimes e^{k\ell m}. \label{Ttensor2}
\eea  
These matrices also obey another $RTT$-type relation: 
\beq
  (\R{q} \otimes 1)\, \uT{E,e}^{(1)}\,  \uT{E,e}^{(2)} 
  =
  \uT{E,e}^{(2)}\,  \uT{E,e}^{(1)}\, (\R{q} \otimes 1).
  \label{RTT2}
\eeq
Application of the tensored map 
$ \Phi \otimes {\rm id} \otimes {\rm id} $ to 
(\ref{RTT1}) converts the $RTT$-type relation into the 
Yang-Baxter equation of the form
\beq
  \R{q\,23}  \R{q\,12}  \R{q\,13} = 
  \R{q\,13}  \R{q\,12}  \R{q\,23}, \label{YB1}
\eeq
while similar action of the conjugate map 
$ {\rm id} \otimes {\rm id} \otimes \Psi $ on (\ref{RTT2}) 
provides the another form of Yang-Baxter equation:
\beq
   \R{q\,12}  \R{q\,13}  \R{q\,23}
   =
   \R{q\,23} \R{q\,13} \R{q\,12}.
   \label{YB2}
\eeq

  Mappings from a universal ${\cal T}$-matrix to a universal $R$-matrix 
has been discussed for only a few quantum algebras. 
Fr\o nsdal\ocite{Fro}  considered such mappings for quantum $gl(n)$ 
and the particular case of two-parametric quantum $ gl(2) $ 
is discussed in Ref.\cite{VJ}. The maps for Alexander-Conway 
quantum algebra is studied in Ref.\cite{CJ2}.

%%%%%%%%%%%%%%%%%%%%%%%%%%%%%%%%%%%%%%%%%%%%%%%%%%
%
%  Representations of A
%
%%%%%%%%%%%%%%%%%%%%%%%%%%%%%%%%%%%%%%%%%%%%%%%%%%
%
\setcounter{equation}{0}
\section{REPRESENTATIONS OF \boldmath{$\A$}}
\label{RepA}

We do not yet have the explicit formulae of finite dimensional 
representation matrices of the function algebra $\A$. These expressions 
will be derived in the next section. But, prior to that, the general 
properties of such finite dimensional representation matrices may be 
understood via the duality arguments interrelating $\A$ and $\U$
algebras. We will address to this task in the present section.

  To be explicit, for us the representations of  $\A,$ signify the 
matrix elements of the universal ${\cal T}$-matrix on $ V^{(\ell)}$ 
defined in \S II:
\bea
  \Tm{\ell}{m'm}(\lambda) &=& (\e{\ell}{m'}(\lambda),\uT{e,E}\, 
  \e{\ell}{m}(\lambda)) \nn \\
  &=& \sum_{abc} (-1)^{(a+c)(a+c-1)/2+(a+c)(\ell-m'+\lambda)} e^{abc}\,
  (\e{\ell}{m'}(\lambda), E_{abc}\, \e{\ell}{m}(\lambda)).
  \label{Tmatdef}
\eea
Under the assumption of the completeness of the basis vectors 
$ \e{\ell}{m}(\lambda), $
it is not difficult to verify the relations:
\beq
  \Delta(\Tm{\ell}{m'm}(\lambda)) = \sum_k \Tm{\ell}{m'k}(\lambda) \otimes 
  \Tm{\ell}{km}(\lambda),
  \qquad
  \epsilon(\Tm{\ell}{m'm}(\lambda)) = \delta_{m'm}.
  \label{corepcond}
\eeq
The relations in (\ref{corepcond}) imply that the matrix elements (\ref{Tmatdef}) 
satisfy the axiom of comodule\ocite{KS}. We may, therefore, regard 
$ \Tm{\ell}{m'm}(\lambda) $ as the 
$ 2\ell+1 $ dimensional matrix representation of the algebra $\A.$ 

  We now consider a product of two representations in order to 
obtain their composition rule. 
We evaluate the matrix elements of $ \uT{e,\Delta(E)} $ on the 
coupled basis vector space in two different ways. 
The first evaluation is as follows:
\bea
   & & (\e{\ell'}{m'}(\ell_1,\ell_2,\Lambda), \uT{e,\Delta(E)}\,
    \e{\ell}{m}(\ell_1,\ell_2,\Lambda))
   \nn \\
   & & \qquad =
   \sum_{abc} (-1)^{(a+c)(a+c-1)/2+(a+c)(\ell'-m'+\Lambda)} e^{abc}\, 
   (\e{\ell'}{m'}(\ell_1,\ell_2,\Lambda), \Delta(E_{abc})\,
    \e{\ell}{m}(\ell_1,\ell_2,\Lambda))
   \nn \\
   & & \qquad =
   \delta_{\ell'\ell} (-1)^{(\ell'-m'+\lambda)(\ell_1+\ell_2+\ell'+\lambda)}\,
   \Tm{\ell}{m'm}(\Lambda).
   \label{PR1}
\eea
In the last equality, the pseudo orthogonality (\ref{Normprod}) 
of the coupled basis has been used. An alternate evaluation of the 
said matrix element explicitly uses the Clebsch-Gordan coupling 
of the basis vectors. With the aid of the relation
\[
  \uT{e,\Delta(E)} = \sum_{prq \atop p'q'r'} 
  (-1)^{(p+r+p'+r')(p+r+p'+r'-1)/2}e^{pqr} e^{p'q'r'} \otimes
  E_{pqr} \otimes E_{p'q'r'},
\]
we obtain
\bea
   & & (\e{\ell'}{m'}(\ell_1,\ell_2,\Lambda), \uT{e,\Delta(E)}\,
    \e{\ell}{m}(\ell_1,\ell_2,\Lambda))
   \nn \\
   & & \qquad =
  \sum_{m_1,m_2 \atop m_1',m_2'} 
  (-1)^{(\ell_1-m_1+\lambda)(\ell_2-m_2'+\lambda)}
  \CGC{\ell_1\ \; \ell_2\ \;\, \ell'}{m_1' \  m_2' \ m'}
  \CGC{\ell_1\ \ell_2\ \ell}{m_1 \; m_2 \; m}
  \Tm{\ell_1}{m_1'm_1}(\lambda)\, \Tm{\ell_2}{m_2'm_2}(\lambda).
  \label{PR2}
\eea
Since the results of the two evaluations have to be identical, 
we obtain the product law  for the quantum supergroup $\A$ 
of the following form:
\bea
  & & \delta_{\ell'\ell} \Tm{\ell}{m'm}(\Lambda) = 
   (-1)^{(\ell'-m'+\lambda)(\ell_1+\ell_2+\ell'+\lambda)}
  \nn \\
  & & \qquad \times
  \sum_{m_1,m_2 \atop m_1',m_2'} 
  (-1)^{(\ell_1-m_1+\lambda)(\ell_2-m_2'+\lambda)}
  \CGC{\ell_1\ \; \ell_2\ \;\, \ell'}{m_1' \  m_2' \ m'}
  \CGC{\ell_1\ \ell_2\ \ell}{m_1 \; m_2 \; m}
  \Tm{\ell_1}{m_1'm_1}(\lambda)\, \Tm{\ell_2}{m_2'm_2}(\lambda).
  \label{productlaw}
\eea
Another derivation of the product law (\ref{productlaw}) is 
found in Ref.\cite{AC05}. 

  Two alternate forms of the product law are readily derived 
by using the pseudo orthogonality of CGC (\ref{CGCortho1}) 
and (\ref{CGCortho2}): 
\bea
  & &   \sum_{m'}\CGC{\ell_1\ \ell_2 \ \ell}{n_1\ n_2\ m'} \Tm{\ell}{m'm}(\Lambda)
  \nn \\
  & & \hspace{2cm} =
  \sum_{m_1,m_2} (-1)^{(n_1+m_1)(\ell_2-n_2+\lambda)}
  \CGC{\ell_1\ \ell_2 \ \ell}{m_1\ m_2\ m} \Tm{\ell_1}{n_1m_1}(\lambda)\,
  \Tm{\ell_2}{n_2m_2}(\lambda),
  \label{PRal1} \\
  & & \sum_m (-1)^{(m'+m)(\ell_1+\ell_2+\ell'+\lambda)} 
  \CGC{\ell_1\ \ell_2\ \ell'}{n_1\ n_2\ m} \Tm{\ell'}{m'm}(\Lambda)
  \nn \\
  & & \hspace{2cm}  =
  \sum_{m_1',m_2'}(-1)^{(n_2+m_2')(\ell_1-n_1+\lambda)}
  \CGC{\ell_1\ \ell_2\ \ell'}{m_1'\; m_2'\; m'}
  \Tm{\ell_1}{m_1'n_1}(\lambda)\,\Tm{\ell_2}{m_2'n_2}(\lambda).
  \label{PRal2}
\eea

  The product law allows us to derive the orthogonality and 
the recurrence relations of the representation matrix $\Tm{\ell}{m'm}. $ 
Setting $ \ell_1 = \ell_2,\ \ell = m = 0 $ in  
(\ref{PRal1}) and using the formula of CGC given in Ref.\cite{AC05},
one can verify  the orthogonality relation:
\beq
  \sum_m (-1)^{m_1(m_1+m) + m_1(m_1-1)/2 + m(m-1)/2} 
  q^{(m_1-m)/2} \Tm{\ell}{m_1m}(\lambda) \Tm{\ell}{-m_2 -m}(\lambda)
  = \delta_{m_1m_2}.
  \label{Tortho0}
\eeq
Another orthogonality relation is similarly obtained by   
setting $ \ell_1 = \ell_2,\ \ell' = m' = 0 $ in the product law 
(\ref{PRal2}):
\beq
 \sum_m (-1)^{(m_1+m)m_1 + m_1(m_1-1)/2 + m(m-1)/2} q^{(m_1-m)/2}\,
 \Tm{\ell}{mm_1}(\lambda)\, \Tm{\ell}{-m-m_2}(\lambda) = \delta_{m_1m_2}.
 \label{Tortho}
\eeq
For the choice of $ \ell_2 = 1 $ in (\ref{PRal1}), the recurrence 
relations for the representation matrices are obtained below. These 
relations are classified into three sets according to the values 
of $ \ell$:\\ 
$\bullet$\,\,The first set has the value of $ \ell = \ell_1 +1.$ It 
comprises of three relations corresponding to all possible values of $n_2.$ 
The recurrence relations listed below correspond to $ n_2 = 1, 0 $ and 
$-1,$ respectively:
\setlength{\jot}{3mm}
\bea
  & & (-1)^{(n+m)\lambda} q^{-(\ell-n)/2} \Fl{n}{0}{-1} 
      \Tm{\ell}{nm}(\Lambda)
      =
      q^{-(\ell-m)/2} \Fl{m}{0}{-1} \Tm{\ell-1}{n-1\,m-1}(\lambda)\, a
  \nn \\
  & & \hspace{2cm}
      -(-1)^{\lambda} q^{m/2} \Gl{m}{0}{0} 
      \Tm{\ell-1}{n-1\,m}(\lambda)\, \alpha
     + q^{(\ell+m)/2} \Fl{-m}{0}{-1} \Tm{\ell-1}{n-1\,m+1}(\lambda)\, b,
  \nn \\
  %%%%%%%%%
  & & (-1)^{(n+m)(1+\lambda)} q^{n/2} \Gl{n}{0}{0}
      \Tm{\ell}{nm}(\Lambda)
      =
      -q^{-(\ell-m)/2} \Fl{m}{0}{-1} \Tm{\ell-1}{n\,m-1}(\lambda)\, \gamma
  \nn \\
  & & \hspace{2cm}
      + q^{m/2} \Gl{m}{0}{0} \Tm{\ell-1}{nm}(\lambda)\, e
     - q^{(\ell+m)/2} \Fl{-m}{0}{-1} \Tm{\ell-1}{n\,m+1}(\lambda)\, \beta,
  \label{recT1} \\
  %%%%%%
  & & (-1)^{(n+m)(1+\lambda)} q^{(\ell+n)/2} \Fl{-n}{0}{-1} 
      \Tm{\ell}{nm}(\Lambda)
      =
      q^{-(\ell-m)/2} \Fl{m}{0}{-1} \Tm{\ell-1}{n+1\,m-1}(\lambda)\, c
  \nn \\
  & & \hspace{2cm}
      + q^{m/2} \Gl{m}{0}{0}
      \Tm{\ell-1}{n+1\,m}(\lambda)\, \delta
      + q^{(\ell+m)/2} \Fl{-m}{0}{-1} \Tm{\ell-1}{n+1\,m+1}(\lambda)\, \beta,
  \nn   
\eea
where
\bea
  & & \Fl{m}{a}{b} = \sqrt{\ku{\ell+m+a}\ku{\ell+m+b}}, \nn \\
  & & \Gl{m}{a}{b} = \sqrt{\ku{2}\ku{\ell+m+a}\ku{\ell-m+b}}, \nn
\eea
and the matrix elements for the fundamental representation $(\ell = 1)$ 
are denoted as 
\beq
   \begin{pmatrix}
     a & \alpha & b \\
     \gamma & e & \beta \\
     c & \delta & d
   \end{pmatrix}
   =
   \left(
   \begin{array}{lll}
     \Tm{1}{11}(\lambda) & \Tm{1}{10}(\lambda) & \Tm{1}{1-1}(\lambda) \\
     \Tm{1}{01}(\lambda) & \Tm{1}{00}(\lambda) & \Tm{1}{0-1}(\lambda) \\
     \Tm{1}{-11}(\lambda) & \Tm{1}{-10}(\lambda) & \Tm{1}{-1-1}(\lambda)
   \end{array}\right).
   \label{T1osp}
\eeq
In this set, the highest weight $ \Lambda  $ assumes a constant value of 
0 (mod $2$). \\
$\bullet$\,\,
The second set corresponds to $ \ell = \ell_1. $ It also contains three 
recurrence relations and, in this instance, we have
$ \Lambda = 1 $ (mod $2$):
\bea
  & & (-1)^{\ell-n+\lambda + (n+m+1)\lambda} q^{(n-m)/2} 
  \Gl{n}{0}{1} \Tm{\ell}{nm}(\Lambda)
  =
  (-1)^{\ell-m} \Gl{m}{0}{1} \Tm{\ell}{n-1\,m-1}(\lambda)\, a
  \nn \\
  & & \hspace{2cm}
  - H^{\ell}_m \, \Tm{\ell}{n-1\,m}(\lambda)\, \alpha
  + \ku{2}^{-1/2} \Gl{m}{1}{0}  \Tm{\ell}{n-1\,m+1}(\lambda)\, b
  \nn \\
  %%%%%%
  & & (-1)^{(n+m)(1+\lambda)} q^{(n-m)/2} H^{\ell}_n\, \Tm{\ell}{nm}(\Lambda)
      = 
      (-1)^{\ell-m} \Gl{m}{0}{1} \Tm{\ell}{n\,m-1}(\lambda)\,\gamma
  \nn \\
  & & \hspace{2cm}
  + H^{\ell}_m\, \Tm{\ell}{nm}(\lambda)\, e 
  + \ku{2}^{-1/2} \Gl{m}{1}{0} \Tm{\ell}{n\,m+1}(\lambda)\, \beta,
  \label{recT2} \\
  %%%%%%%
  & & (-1)^{(n+m)\lambda} q^{(n-m)/2} \ku{2}^{-1/2} \Gl{n}{0}{1} \Tm{\ell}{nm}(\Lambda)
  = 
  (-1)^{\ell-m} \Gl{m}{0}{1} \Tm{\ell}{n+1\,m-1}(\lambda)\, c
  \nn \\
  & & \hspace{2cm}
  - H^{\ell}_m \, \Tm{\ell}{n+1\,m}(\lambda)\, \delta
  + \ku{2}^{-1/2} \Gl{m}{1}{0} \Tm{\ell}{n+1\,m+1}(\lambda)\, d,
  \nn
\eea
where
\[
  H^{\ell}_m = q^{-\ell/2} \ku{\ell+m+1} - (-1)^{\ell-m} q^{\ell/2} \ku{\ell-m+1}.
\]
$\bullet$\,\, Similarly, the third set contains three recurrence 
relations  for $ \ell = \ell_1 - 1. $ For this case,
the highest weight reads $ \Lambda = 0 $ (mod $2$), and the 
recurrence relations are given by
\bea
   & & q^{(\ell-m+n+1)/2} \Fl{-n}{1}{2} \Tm{\ell}{nm}(\Lambda)
   =
   q^{(\ell+1)/2} \Fl{-m}{1}{2} \Tm{\ell+1}{n-1\,m-1}(\lambda)\, a
   \nn \\
   & & \hspace{2cm}
   + (-1)^{\ell-m+\lambda} \Gl{m}{1}{1} \Tm{\ell+1}{n-1\,m}(\lambda)\, \alpha
   - q^{-(\ell+1)/2} \Fl{m}{1}{2} \Tm{\ell+1}{n-1\,m}(\lambda)\, b,
   \nn \\
   %%%%%%%%%%
   & & (-1)^{\ell-n+\lambda} q^{(n-m)/2} \Gl{n}{1}{1} \Tm{\ell}{nm}(\Lambda)
   =
   q^{(\ell+1)/2} \Fl{-m}{1}{2} \Tm{\ell+1}{n\,m-1}(\lambda)\, \gamma
   \label{recT3} \\
   & & \hspace{2cm}
   + (-1)^{\ell-m+\lambda} \Gl{m}{1}{1} \Tm{\ell+1}{nm}(\lambda)\, e 
   - q^{-(\ell+1)/2} \Fl{m}{1}{2} \Tm{\ell+1}{n\,m+1}(\lambda)\, \beta,
   \nn \\ 
   %%%%%%%%%%%%
   & & q^{-(\ell-n+m+1)/2} \Fl{n}{1}{2} \Tm{\ell}{nm}(\Lambda)
   =
   -q^{(\ell+1)/2} \Fl{-m}{1}{2} \Tm{\ell+1}{n+1\,m-1}(\lambda)\, c
   \nn \\
   & & \hspace{2cm}
   -(-1)^{\ell-m+\lambda} \Gl{m}{1}{1} \Tm{\ell+1}{n+1\,m}(\lambda)\, \delta
   + q^{-(\ell+1)/2} \Fl{m}{1}{2} \Tm{\ell+1}{n+1\,m+1}(\lambda)\, d.
   \nn
\eea

  The discussion so far is independent of the explicit form of the 
universal ${\cal T}$-matrix. The general properties of the universal 
${\cal T}$-matrix and the the Clebsch-Gordan decomposition 
of tensor product representations play a seminal role 
in the derivation of all properties of the representation of $\A.$  
Thus one can repeat the same arguments for other quantum deformations of 
$ OSp(1/2),$ namely the Jordanian\ocite{CK98} and the 
super-Jordanian\ocite{ACS2,BLT} analogs for
deriving their product law, the orthogonality and the recurrence 
relations. Being triangular algebras, the Jordanian and the 
super-Jordanian deformations of $OSp(1/2)$ possess the 
same Clebsch-Gordan decomposition as in the present case. 

  The representations of $ SL_q(2) $ (or $ SU_q(2)$) have been 
discussed by many authors. Among the properties of representation matrices, 
the product law\ocite{GKK}, recurrence relations\ocite{GKK,Nom}, orthogonality 
and $RTT$-relation\ocite{Nom} and generating functions\ocite{Nom2} 
are found in literature. In Ref.\cite{Nom}, the representation matrices 
are interpreted as the wave functions of quantum symmetric top in 
noncommutative space. Representations of the  Jordanian 
quantum group $SL_{\sf h}(2)$ have been considered in 
Refs.\cite{CQ,NA00}.

%%%%%%%%%%%%%%%%%%%%%%%%%%%%%%%%%%%%%%%%%%%%%%%%
%
%  Representation matrix and polynomial
%
%%%%%%%%%%%%%%%%%%%%%%%%%%%%%%%%%%%%%%%%%%%%%%%
%
\setcounter{equation}{0}
\section{REPRESENTATION MATRIX AND LITTLE \boldmath{$Q$}-JACOBI POLYNOMIALS}
\label{RepPol}

  Explicit formulae for the representation matrices of the quantum group  
$SU_q(2)$ have been obtained by several authors\ocite{VS,Koo,MMNNU}. 
It is observed that, for the finite dimensional representations, 
the matrix elements are expressed in terms of the little $q$-Jacobi 
polynomials. Investigating the Jordanian quantum group $SL_{\sf h}(2)$ 
in a similar framework, it has also been noted that the conventional 
Jacobi polynomials contribute\ocite{NA00} to the representation 
matrices therein. The corresponding matrix elements for the two-parametric 
quantum group $GL(2)$ have been computed in Ref.\cite{JV}. Their 
relation to orthogonal polynomials, however, is still an open problem.   

 In this section, we provide the explicit form of the representation 
matrices (\ref{Tmatdef}) of $\A$ by direct computation, and study the 
resulting polynomial structure. Towards this end, we proceed by noticing 
the identities obtained by repeated use of (\ref{RepU}):
\bea
  & & V_-^c\, \e{\ell}{m}(\lambda) = (-1)^{c(\ell-m) + c(c-3)/2}\left(
   \frac{1}{\ku{2}^c}
   \frac{\ku{\ell+m}!}{\ku{\ell-m}!} \frac{\ku{\ell-m+c}!}{\ku{\ell+m-c}!}
  \right)^{1/2} \e{\ell}{m-c}(\lambda),
  \nn \\
  & & V_+^a\, \e{\ell}{m}(\lambda) = \left(
   \frac{1}{\ku{2}^a}
   \frac{\ku{\ell-m}!}{\ku{\ell+m}!} \frac{\ku{\ell+m+a}!}{\ku{\ell-m-a}!}
  \right)^{1/2} \e{\ell}{m+a}(\lambda). \label{vpe}
\eea
The explicit listing of the basis elements of $\A$ in (\ref{enrs}) 
renders the computation of the matrix elements straightforward. The  
final result is quoted below:
\bea
  & & \Tm{\ell}{m'm}(\lambda) \nn \\
  & & = 
  (-1)^{(m'-m)(m'-m-1)/2 + (m'-m)(\ell-m'+\lambda)}
  \frac{q^{m(m'-m)/2}}{\sqrt{\ku{2}^{m'-m}}}
  \left(
     \frac{\ku{\ell+m}!}{\ku{\ell-m}!}
     \frac{\ku{\ell+m'}!}{\ku{\ell-m'}!}
  \right)^{1/2}
  \nn \\ 
  & &  \times
  \sum_c (-1)^{c(\ell-m)} \frac{q^{-c(m'-m)/2}}{\ku{2}^c}
  \frac{\ku{\ell-m+c}!}{\ku{\ell+m-c}!}
  \frac{x^{m'-m+c}}{\ku{m'-m+c}!} 
  \exp\left( \frac{m-c}{2}z \right)
  \frac{y^c}{\ku{c}!},
  \nn \\
  & & \hfill \label{RepMatA}
\eea
where the index $c$ runs over all non-negative integers maintaining the 
argument of $\{{\cal X}\}_q$ non-negative. 

  The fundamental representation ($\ell = 1, \ \lambda = 0$) may be 
identified with the quantum supermatrix of Ref.\cite{AC05} 
{\it \`{a} la} (\ref{T1osp}). 
This identification allows us to realize the entries of 
the quantum supermatrix in terms of the generators of  
the quantum supergroup $\A$:
\bea
  & & a = xy + e^{z/2} +  \frac{x^2 e^{-z/2} y^2}{\ku{2}^2},
  \quad
  \alpha = x - \frac{x^2 e^{-z/2} y}{q^{1/2} \ku{2}},
  \quad
  b = - \frac{x^2 e^{-z/2}}{q\ku{2}},
  \nn \\
  & & \gamma = y + \frac{q^{1/2}x e^{-z/2} y^2}{\ku{2}},
   \quad
   e = 1 - x e^{-z/2} y, \quad \beta = -q^{-1/2} xe^{-z/2},
  \label{OSPtoxyz} \\
  & & c = - \frac{q e^{-z/2}y^2}{\ku{2}}, \qquad\qquad
      \delta = q^{1/2} e^{-z/2} y, \qquad\quad
      d = e^{-z/2}. \nn
\eea
Straightforward computation using the commutation relations 
(\ref{Acomm}) allows us to infer that the realization (\ref{OSPtoxyz}) 
recovers all the commutation relations of the supermatrix listed in 
Ref.\cite{AC05}. The realization (\ref{OSPtoxyz}), more 
importantly, implies the following Gaussian decomposition of the 
quantum supermatrix\ocite{DKLS,DKS,CK96}:
\bea
  & &  
  \begin{pmatrix}
    1 & 0 & 0 \\
    -q^{-1/2}x & 1 & 0 \\
    -\frac{1}{q\ku{2}}x^2 & x & 1
  \end{pmatrix}
  \begin{pmatrix}
    e^{-z/2} & 0 & 0 \\
    0 & 1 & 0 \\
    0 & 0 & e^{z/2}
  \end{pmatrix}
  \begin{pmatrix}
    1 & q^{1/2} y & -\frac{q}{\ku{2}}y^2 \\
    0 & 1 & y \\
    0 & 0 & 1
  \end{pmatrix}
  \nn \\
  & &  =
  \begin{pmatrix}
    d & \delta & c \\ 
    \beta & e & \gamma \\
    b & \alpha & a
  \end{pmatrix} 
  = 
  C \begin{pmatrix}
     a & \alpha & b \\
     \gamma & e & \beta \\
     c & \delta & d
  \end{pmatrix} C^{-1},
  \label{Gdecomp}
\eea
where
\[
  C = 
  \begin{pmatrix}
    0 & 0 & 1 \\ 0 & 1 & 0 \\ 1 & 0 & 0
  \end{pmatrix},
  \qquad
  C^{-1} = C.
\]

   We now turn our attention to the polynomial structure built into the 
general matrix element (\ref{RepMatA}) in terms of the variable 
\beq
  \zeta = \frac{q^{-1/2}}{\ku{2}} xe^{-z/2}y.  \label{zetadef}
\eeq
To demonstrate this, the product of generators in (\ref{RepMatA}) for 
the case $m'-m \geq 0$ may be rearranged as follows:
\[
  x^{m'-m+c} \exp\left( \frac{m-c}{2}z \right) y^c 
  =
  (-1)^{c(c-1)/2} q^{-mc} x^{m'-m} e^{mz/2} (x e^{-z/2} y)^c.
\]
The matrix element $\Tm{\ell}{m'm}(\lambda)$ may now be succinctly
expressed as a polynomial structure given below: 
\bea
  \Tm{\ell}{m'm}(\lambda) &= & 
  (-1)^{(m'-m)(m'-m-1)/2 + (m'-m)(\ell-m'+\lambda)}
  \frac{q^{m(m'-m)/2 }}{\ku{m'-m}!\sqrt{\ku{2}^{m'-m}}}
  \nn \\
  & \times &
  \left(
     \frac{\ku{\ell-m}!}{\ku{\ell+m}!}
     \frac{\ku{\ell+m'}!}{\ku{\ell-m'}!}
  \right)^{1/2}
  x^{m'-m} e^{mz/2} \Pz{\ell}{m'm}.
  \label{Tmatposi}
\eea
The polynomial $ \Pz{\ell}{m'm}$ in the variable $\zeta$ is defined by
\bea
  & & \Pz{\ell}{m'm} =
  \sum_c (-1)^{c(\ell-m)+c(c-1)/2} q^{-c(m'+m-1)/2}
  \nn \\
  & & \hspace{2cm} \times
  \frac{\ku{m'-m}! \ku{l+m}! \ku{\ell-m+c}!}
   {\ku{m'-m+c}! \ku{\ell+m-c}! \ku{\ell-m}! \ku{c}!}\,
  \zeta^c,
  \label{Pzposi}
\eea
where the index $c$ runs over all non-negative integers maintaining the 
arguments of $\{{\cal X}\}_q$ non-negative. For the case 
$ m'-m \leq 0,$ we make a replacement of the summation index $c$ with  
$a = m'-m+c$. Rearrangement of the generators now provides the 
following expression of the general matrix element:
\bea
  \Tm{\ell}{m'm}(\lambda) &= & 
  (-1)^{(m-m')(m-m'+1)/2 - (m-m')(\lambda-1)}
  \frac{q^{-m'(m-m')/2}}{\ku{m-m'}!\sqrt{\ku{2}^{m-m'}}}
  \nn \\
  & \times &
  \left(
     \frac{\ku{\ell+m}!}{\ku{\ell-m}!}
     \frac{\ku{\ell-m'}!}{\ku{\ell+m'}!}
  \right)^{1/2} e^{m'z/2}\,y^{m - m'}\,
   \Pz{\ell}{m'm},
  \label{Tmatnega}
\eea
where the polynomial $ \Pz{\ell}{m'm} $ for $ m'-m \leq 0 $ is defined by
\bea
  & & \Pz{\ell}{m'm} = 
  \sum_a (-1)^{a(\ell-m')+a(a-1)/2} q^{-a(m'+m-1)/2}
  \nn \\
  & & \hspace{2cm} \times \frac{\ku{m-m'}! \ku{l+m'}! \ku{\ell-m'+a}!}
   {\ku{m-m'+a}! \ku{\ell+m'-a}! \ku{\ell-m'}! \ku{a}!}\,
  \zeta^a.
  \label{Pznega}
\eea
It is immediate to note that the polynomials are symmetric 
with respect to the transposition $m \leftrightarrow m^{\prime}$:  
$
  \Pz{\ell}{mm'} = \Pz{\ell}{m'm},
$ 
and that 
$ \Pz{\ell}{-\ell\;m} = \Pz{\ell}{m'\;-\ell} = 1. $ 

\par

Polynomials obtained above are related to the basic hypergeometric functions. 
We define the basic hypergeometric function $ {}_2 \phi_1 $ by\ocite{GR}
\beq
  {}_2 \phi_1(a_1,a_2; b; Q; z) 
  = 
  \sum_{n=0} 
  \frac{(a_1;Q)_n (a_2;Q)_n}{(b;Q)_n (Q;Q)_n}z^n,
  \label{F21}
\eeq
where the shifted factorial is defined as usual:
\beq
  (x;Q)_n = \left\{
    \begin{array}{cl}
        1, & \quad n=0 \\ 
        {\displaystyle \prod_{k=0}^{n-1} (1-xQ^k)}. & \quad n \neq 0
    \end{array}
  \right.
  \label{siftf}
\eeq 
The little Q-Jacobi polynomials are defined via $ {}_2 \phi_1 $ as 
standard theory of orthogonal polynomials\ocite{GR}
\beq
  p^{(a,b)}_m(z) = {}_2 \phi_1(Q^{-m}, abQ^{m+1}; aQ; Q; Qz).
  \label{DefqJ}
\eeq
Setting $ a = Q^{\alpha}, b = Q^{\beta}, $ 
we have the following form of little $Q$-Jacobi polynomials
\beq
  p^{(\alpha,\beta)}_m(z) = 
  \sum_{n}  
  \frac{(Q^{-m};Q)_n (Q^{\alpha+\beta+m+1};Q)_n}{(Q^{\alpha+1};Q)_n (Q;Q)_n}
  (Qz)^n.
  \label{DefqJ2}
\eeq
Rewriting our polynomials (\ref{Pzposi}) and (\ref{Pznega}) 
in terms of the shifted factorial with $ Q = -q, $ one can 
identify our polynomials with the little Q-Jacobi polynomials. 
For the choice $ m'-m \geq 0,$ the polynomial structure reads 
\beq
  P^{\ell}_{m'm}(\zeta) =
     \sum_a  \,
     \frac{((-q)^{-\ell-m};-q)_a\, ((-q)^{\ell-m+1};-q)_a}
          {((-q)^{m'-m+1};-q)_a\, (-q;-q)_a} \,(-q\zeta)^a
     = p^{(m'-m,-m'-m)}_{\ell+m}(\zeta),
  \label{PtoJ1}
\eeq
and for the $ m'-m \leq 0$ case its identification is given by  
\beq
  \displaystyle
  P^{\ell}_{m'm}(\zeta) =
     \sum_a \,
     \frac{((-q)^{-\ell-m'};-q)_a \, ((-q)^{\ell-m'+1};-q)_a}
          {((-q)^{m-m'+1};-q)_a\, (-q;-q)_a} \,(-q\zeta)^a
     = p^{(m-m',-m'-m)}_{\ell+m'}(\zeta).
  \label{PtoJ2}
\eeq
It is amazing that $Q=-q$ polynomials appear for the supergroup $\A$ 
in contrast to the $Q=q$ polynomials being present for the quantum  
group $SU_q(2)$\ocite{VS,Koo,MMNNU}. 

%%%%%%%%%%%%%%%%%%%%%%%%%%%%%%%%%%%%%%%%%%%%%%%%%%
%
%  Concluding Remarks
%
%%%%%%%%%%%%%%%%%%%%%%%%%%%%%%%%%%%%%%%%%%%%%%%%%%
%
\section{CONCLUDING REMARKS}
\label{ConclR}

  Starting from the construction of the universal ${\cal T}$-matrix, 
we have investigated the finite dimensional representations of
the quantum supergroup $\A.$ A qualitative difference between 
the universal ${\cal T}$-matrices for the classical and the quantum 
$OSp(1/2)$ algebras exists due to the nilpotency of the classical parity 
odd elements. The absence of the said nilpotency in the quantum 
case induces a new polynomial structure in the matrix elements of 
${\cal T}$. We observe that these polynomials are expressed in 
terms of the little $Q$-Jacobi polynomials.
This suggests a new link between orthogonal polynomials 
and representations of quantum supergroups. It is likely to be a 
general property that if the Grassmann variables in a classical 
supergroup lose nilpotency at the quantum level, it may be reflected 
in the representations of the quantum supergroup so that the entries 
of representation matrices may have a new quantized polynomial 
structure. The present work provides an example of this statement. 
Another likely candidate for the existence of similar polynomial 
structure is the super-Jordanian $OSp(1/2)$ where 
the loss of nilpotency has been observed\ocite{AC04}.
We believe that the investigation along this line will 
give a new algebraic background to basic hypergeometric series. 

  We have also tried to extend the known properties of 
the representations for $ SL_q(2)$ to the quantum supergroup $\A.$ 
An extension of the product law, the orthogonality and the recurrence 
relations was shown to be possible. 
However, two known results are not extended  in the present 
work, that is, the generating function\ocite{Nom2} of the 
representation matrices and the Peter-Weyl theorem\ocite{Wor,MMNNU}. 
In order to discuss the Peter-Weyl theorem, a Haar measure 
has to be defined on $\A$. Since the loss of nilpotency 
makes the superspace on $\A$ more complex than the classical case, 
studying the Peter-Weyl theorem may be interesting from the 
viewpoint of harmonic analysis on quantum supergroups. 
We will discuss these issues elsewhere. 

%%%%%%%%%%%%%%%%%%%%%%%%%%%%%%%%%%%%%%%%%%%%%%%%%%%
%
%    Acknowledgements
%
%%%%%%%%%%%%%%%%%%%%%%%%%%%%%%%%%%%%%%%%%%%%%%%%%%%
%
\section*{ACKNOWLEDGEMENTS}
A comment by the Referee motivated us in finding the present form of  
(\ref{uniTclassical}). We thank him for this.  
One of us (N.A.) would like to thank R. Jagannathan 
for his warm hospitality at The Institute of Mathematical 
Sciences where a part of this work was done. 
The work of N.A. is partially supported by 
the grants-in-aid from JSPS, Japan (Contract No. 15540132).
Other authors (R.C., S.S.N.M. and J.S.) are partially supported by
the grant DAE/2001/37/12/BRNS, Government of India.

%%%%%%%%%%%%%%%%%%%%%%%%%%%%%%%%%%%%%%%%%%%%%%%%%%
%
%   Appendix
%
%%%%%%%%%%%%%%%%%%%%%%%%%%%%%%%%%%%%%%%%%%%%%%%%%%%
%
\setcounter{equation}{0}
\section*{APPENDIX: UNIVERSAL ${\cal T}$-MATRIX FOR QUANTUM SUPERGROUPS}
\renewcommand{\theequation}{A.\arabic{equation}}

This Appendix is devoted to a general discussion of the universal 
${\cal T}$-matrix for quantum supergroups. In particular, the 
relations used in \S III and \S IV are proved in general setting. 

  Let $ \U $ and $ \A $ be dually conjugate unital $ {\mathbb Z}_2 $ 
graded Hopf algebras. 
The basis of the algebras $ \U $ and $ \A $ are denoted by 
$ E_k, e^k, $ respectively. One may assume that 
$ E_0 = 1_{\U},\ e^0 = 1_{\A} $  and $ p(E_k) = p(e^k) $ 
without loss of generality. 
The duality of $ \U $ and $ \A $ is reflected in the structure 
constants:
\bea
  & & E_k E_{\ell} = \sum_m f_{k \ \ell}^{\; m} \, E_m, \qquad
  \Delta(E_k) = \sum_{pq} g^{p\ q}_{\;k} E_p \otimes E_q,
  \label{EkEl} \\
  & & e^k e^{\ell} = \sum_m g^{k\ \ell}_{\;m} \, e^m, \qquad\quad 
  \Delta(e^k) = \sum_{pq} f_{p\ q}^{\;k} e^p \otimes e^q.
  \label{ekel}
\eea
The universal ${\cal T}$-matrix is defined by
\beq
  \uT{e,E} = \sum_k (-1)^{p(e^k)(p(e^k)-1)/2} e^k \otimes E_k, \quad
  \in {\cal A} \otimes {\cal U}
  \label{uTdef}
\eeq
Although the factor $ (-1)^{p(e^k)(p(e^k)-1)/2} $ is trivial, it is 
convenient to keep it in the discussion of universal ${\cal T}$-matrix. 

  We start with  the proof of the relations in (\ref{uniTprop}). 
The proof of the first relation in (\ref{uniTprop}) is straightforward
\bea
  \uT{e,E} \uT{e',E} 
  &=& \sum_{k,\ell,m} (-1)^{(p(e^k)+p(e^{\ell})) (p(e^k)+p(e^{\ell})-1)/2} e^k \otimes e^{\ell} 
      \otimes f_{k\; \ell}^{\;m} E_m
  \nn \\
  &=& \sum_{m} (-1)^{p(e^m)(p(e^m)-1)/2} \Delta(e^m) \otimes E_m
  \nn \\
  &=& \uT{\Delta(e),E}. \nn
\eea
The second equality is due to the fact that the parity of the both sides 
of the first equation in (\ref{EkEl}) are equal, 
and that the structure constants are of even parity. 
The second relation in (\ref{uniTprop}) can be proved similarly. The third 
relation in (\ref{uniTprop}) follows from 
$ \epsilon(e^k) = \delta^k_0,\; \epsilon(E_k) = \delta_{k0}. $ 
The last relation in (\ref{uniTprop}) is a consequence of the identities:
\beq
  \uT{e,S(E)} \uT{e,E} = \uT{e,E} \uT{e,S(E)} = 1 \otimes 1, \quad
  \uT{S(e),E} \uT{e,E} = \uT{e,E} \uT{S(e),E} = 1 \otimes 1. \label{Sinverse}
\eeq
The above identities can be proved by using the axiom of antipode and 
$ p(S(e^k)) = p(e^k).$ 

   We next derive the $RTT$-type relations (\ref{RTT1}) and 
(\ref{RTT2}). 
Defining
\[
 \uT{e,E}^{(1)} = \sum_k (-1)^{p(e^k)(p(e^k)-1)/2} e^k \otimes E_k \otimes 1,
 \quad
 \uT{e,E}^{(2)} = \sum_k (-1)^{p(e^k)(p(e^k)-1)/2} e^k \otimes 1 \otimes E_k,
\]
we obtain
\[
 \uT{e,E}^{(1)}\, \uT{e,E}^{(2)} 
 = \sum_m (-1)^{p(e^m)(p(e^m)-1)/2} e^m \otimes \Delta(E_m).
\]
On the other hand the transposed coproduct appears in the 
reverse-ordered product:
\bea
  \uT{e,E}^{(2)}\, \uT{e,E}^{(1)} 
  &=& \sum_m (-1)^{p(e^m)(p(e^m)-1)/2} e^m \otimes \Delta'(E_m)
  \nn \\
  &=& \sum_m (-1)^{p(e^m)(p(e^m)-1)/2} e^m \otimes {\cal R} \Delta(E_m) {\cal R}^{-1} 
  \nn \\
  &=& (1 \otimes {\cal R})\, \uT{e,E}^{(1)}\, \uT{e,E}^{(2)}\, (1 \otimes {\cal R}^{-1}).
  \nn
\eea
In the last equality, the fact that the universal ${\cal R}$-matrix 
is of even parity has been used. This completes our proof of 
the $RTT$-type relation (\ref{RTT1}). The proof of (\ref{RTT2}) may be 
done similarly with the definitions
\[
 \uT{E,e}^{(1)} = \sum_k (-1)^{p(e^k)(p(e^k)-1)/2} E_k \otimes 1 \otimes e^k,
 \quad
 \uT{E,e}^{(2)} = \sum_k (-1)^{p(e^k)(p(e^k)-1)/2} 1 \otimes  E_k \otimes e^k.
\]

\medskip
\noindent 
{\bf Note added:} After submitting the manuscript, we have learnt about 
the work of Zou \ocite{Z03}, where the representations of $OSp_{q}(1/2)$
and their relation to the basic hypergeometric functions are  
computed by using other basis states, and adopting a method (similar
to Ref. \cite{MMNNU}) different from ours. Peter-Weyl theorem has also 
been established there.
}

%%%%%%%%%%%%%%%%%%%%%%%%%%%%%%%%%%%%%%%%%%%%%%%%%%%
%
%    References
%
%%%%%%%%%%%%%%%%%%%%%%%%%%%%%%%%%%%%%%%%%%%%%%%%%%%
%

\end{document}